\documentclass[openacc]{rsproca_new}

\usepackage{lipsum}
\usepackage{amsfonts}
\usepackage{graphicx}
\usepackage{epstopdf}
\usepackage{algorithmic}
\usepackage{bm}
\usepackage{multirow}

\ifpdf
  \DeclareGraphicsExtensions{.eps,.pdf,.png,.jpg}
\else
  \DeclareGraphicsExtensions{.eps}
\fi

\newcommand{\z}{\bm{z}}

\newtheorem{claim}{Claim}
\newtheorem{assumption}{Assumption}
\newtheorem{definition}{Definition}
\newtheorem{corollary}{Corollary}
\newtheorem{theorem}{Theorem}
\newtheorem{lemma}{Lemma}
\newtheorem{proposition}{Proposition}

\begin{document}

\title{
GFINNs: GENERIC Formalism Informed Neural Networks for Deterministic and Stochastic Dynamical Systems}


\author{
Zhen Zhang$^{1}$, Yeonjong Shin$^{1}$ 
and \\
George Em Karniadakis$^{1,2}$}

\address{$^{1}$ Division of Applied Mathematics, 
and $^{2}$ School of Engineering,
Brown University, Providence, RI, 02912, USA\\}

\subject{deep learning, applied mathematics, thermodynamics}


\corres{Yeonjong Shin\\
\email{yeonjong\_shin@brown.edu}}

\begin{abstract}
  We propose the GENERIC formalism informed neural networks (GFINNs)
  that obey the symmetric degeneracy conditions of the GENERIC
  formalism.
  GFINNs comprise two modules, each of which contains
  two components. 
  We model each component using a neural network whose architecture is designed to satisfy the required conditions.
  The component-wise architecture design 
  provides flexible ways of leveraging available physics
  information into neural networks.
  We prove theoretically that GFINNs
  are sufficiently expressive to learn the
  underlying equations, hence establishing the
  universal approxima--
  tion theorem. 
  We demonstrate the performance of GFINNs in three simulation
  problems: gas containers exchanging heat and volume,
  thermoelastic double pendulum and the Langevin 
  dynamics. In all the examples, 
  GFINNs outperform existing methods, hence demonstrating good accuracy in predictions for both deterministic and stochastic systems.
\end{abstract}

\begin{fmtext}
\section{Introduction}
The discovery of governing equations for dynamical systems 
from observed data is a longstanding scientific endeavor \cite{bongard2007automated,brunton2016discovering,schmidt2009distilling}. 
The so-called data-driven discovery dated back to Kepler
refers to scientific methods that extract important features from data 
and either approximate or identify governing equations by means of (parameterized) function classes.
With the recent advancement in deep learning,
neural network classes have been popularly employed in modelling and simulations
and have demonstrated some promising empi--
\end{fmtext}
\maketitle
\noindent 
empirical results 
\cite{chang2017reversible, weinan2017proposal, haber2017stable, lu2018beyond, qin2019data, raissi2018multistep}.

\begin{figure}[ht]
    \centering
    \includegraphics[width=0.9\textwidth]{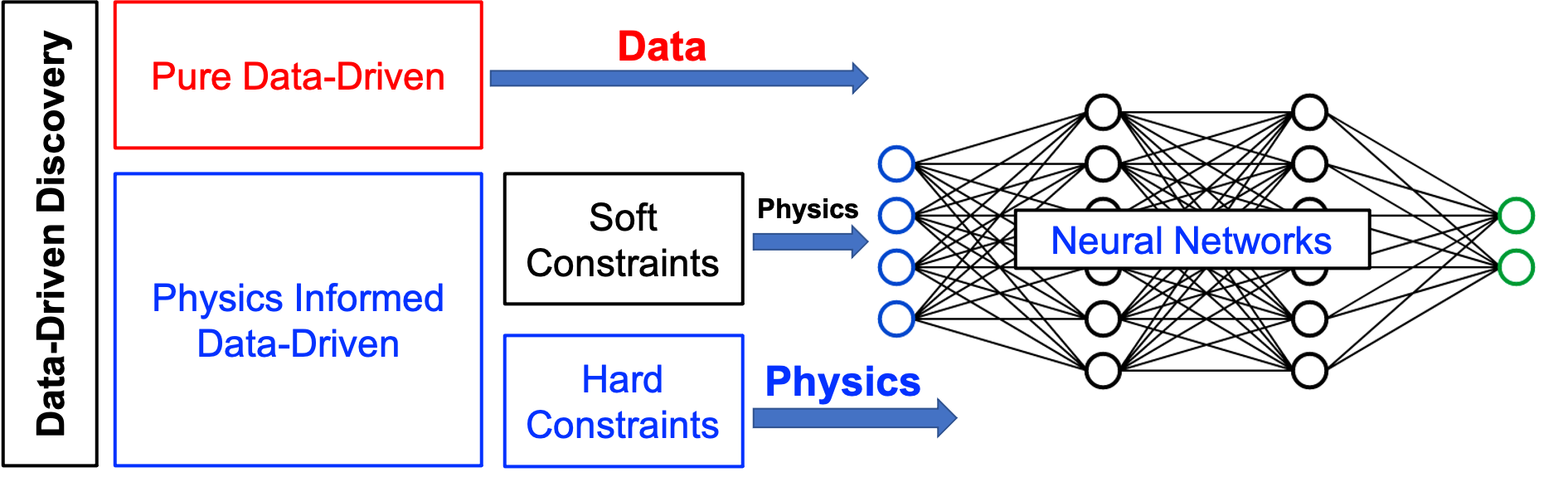}
    \caption{A schematic diagram showing the
    major three approaches in the data-driven discovery of dynamical systems.
    One is the pure data-driven approach.
    The other two are the physics informed data-driven approach.
    The second one is based on soft constraints imposed in an associated loss function.
    The third one directly imposes the underlying physics into neural networks by designing novel network architectures.}
    \label{fig:DD}
\end{figure}

The data-driven discovery may be classified into two major approaches.
One is the pure data-driven methods\cite{ChenRBD18,KidgerMFL20} which model governing equations
using neural networks that are trained to
fit observed data \textit{without physics}. 
This approach could provide a neural network 
model that mimics training trajectories,
particularly 
when no physics but a large amount of data are available.
However, it is likely that the learned models 
do not generalize well on the region where 
training data are scarce or even do not exist.
The other is the physics-informed data-driven approach,
which aims to model governing equations
by \textit{embedding principles of physics
into neural networks} together with data.
By exploiting physics,
it was empirically observed that 
the amount of data needed to get good performance
is much less than those of the pure data-driven methods,
and the learned model is stable and generalizes well.
Typically, the physics are imposed
on neural networks
by means of either soft or hard constraints.
The use of soft constraints 
introduces regularization terms 
in an associated loss function 
that penalize and generalize neural networks
that do not obey the physics \cite{greydanus2019hamiltonian,raissi2019physics,cranmer2020lagrangian,hernandez2021structure}.
The hope is that neural networks 
approximately follow the physics 
after training.
The hard constraints are typically imposed 
by designing proper neural network
architectures that  
obey the underlying principles
without optimization processes \cite{jin2020learning, jin2020sympnets},
yet maintain sufficient expressivity 
so that 
governing equations 
can be learned from data.
This is the approach we follow in the current paper.
A schematic diagram
of the classification of the data-driven discovery
is given in Figure~\ref{fig:DD}.

Among many principles of physics,
we consider the General Equation for Non-Equilibrium Reversible-Irreversible Coupling (GENERIC) formalism \cite{grmela1997dynamics,ottinger1997dynamics,ottinger2005beyond}.
GENERIC provides a general mathematical framework 
describing states beyond equilibrium of a dynamical system \cite{ottinger2005beyond},
which involves two separate generators for the reversible and irreversible dynamics.
These generators are required to satisfy some symmetry and degeneracy conditions,
which constitute a key feature of the GENERIC structure.
These conditions are often interpreted as 
the ﬁrst and second principles of thermodynamics,
which further can be expressed in the language of linear algebra. 


Our goal is to embed the GENERIC structure directly into neural networks, yet to maintain sufficient expressivity.
By leveraging the level of prior physics information under the GENERIC framework,
we propose a systematical approach in designing 
neural network modules.

In the case where either one or all generators are known \textit{a priori},
we design neural network models for either energy or entropy or both
by exploiting certain properties of the generators.
Due to the multiplicative structure of the gradient of neural networks,
care needs to be taken into the input layer
in order to meet the degeneracy conditions.
We thus introduce a transformation in the first layer, which roughly speaking, projects input into a proper low-dimensional space 
on which the degeneracy conditions hold.
On the other hand,
if no prior information is available,
we model generators by
neural networks
whose architectures
are motivated 
by the spectral decomposition of matrix \cite{trefethen1997numerical}.
In any cases,
the proposed neural network models 
obey the required constraints exactly.
We refer to our neural network models 
as the GENERIC formalism informed neural networks (GFINNs).

Furthermore, we prove 
the universal approximation theorem for GFINNs under some assumptions.
Altogether, 
GFINNs not only obey the required GENERIC conditions, but also 
are sufficiently expressive in learning the underlying physical quantities and generators.
We demonstrate the performance of GFINNs 
on several tasks including double pendulums,
gas containers and stochastic differential equations.
We found that GFINNs 
outperform the existing network architectures 
\cite{dietrich2021learning, hernandez2021structure, lee2021machine}
in all the tests we considered.
    
There are two existing works that have attempted to
embed the GENERIC formalism into neural networks.
\cite{hernandez2021structure}
proposed structure preserving 
neural networks (SPNNs),
which aim to learn the physical quantities 
assuming both generators are known.
The degeneracy requirements
are softly constrained in the 
loss function,
which may cause violation 
of the degeneracy conditions 
even after training.
\cite{lee2021machine}
proposed 
the GENERIC neural ordinary differential
equations (GNODEs),
which satisfy the required conditions 
by suitable parameterization of bracket structure,
however,
no universal approximation theorem has been proven. 

The rest of the paper is organized as follows.
Upon introducing the problem setup and 
some preliminaries in Section~\ref{sec:setup},
the GENERIC formalism informed neural networks (GFINNs) are presented in Section~\ref{sec:methods}
along with 
the universal approximation theorem.
Numerical examples are provided 
in Section~\ref{sec:experiments},
demonstrating 
the effectiveness of GFINNs
against other methods.

\section{Problem Setup} \label{sec:setup}
We consider the problem of learning an autonomous dynamical system that can be expressed in the following form:
\begin{equation}
    \dot{\z}(t) = F(\z), \qquad \z \in \Omega \subset \mathbb{R}^d, 
    \quad t \in (0,T], 
    \quad \z(0) = \z_0,
\end{equation}
where $t$ refers to the time coordinates
and 
$F(\z)$ is the unknown vector-valued field function.

With the goal of approximating $F(\z)$ through data,
we employ a neural network $F_{\text{NN}}(\z;\theta)$ 
whose parameters are determined so that 
the trajectories generated from the resulting dynamics $\dot{\z}(t) = F_{\text{NN}}(\z;\theta)$
are close to the observed trajectory data.
More precisely, let $N_{\text{traj}}$ be the number of observed trajectories.
Let $T$ be the number of time-steps 
and $t_j$ be the $j$-th time stamp,
which we set them equal for all trajectories for the sake of notational simplicity.
Let $\z(t;\z_0)$ be the state variable at time $t$ whose value at $t_0$ is $\z_0$.
Then, the $k$-th trajectory is written as $\{t_j, \z(t_j;\z_0^{(k)})\}_{j=0}^{T}$.
We then seek to find the optimal network parameters 
that minimize the loss function defined by
\begin{equation}\label{eq:loss}
    \mathcal{L}(\theta) = \frac{1}{N_{\text{traj}}} \sum_{k=1}^{N_{\text{traj}}} \frac{1}{T} \sum_{j=1}^{T} 
    \left\Vert\z_{\text{NN}}^{\theta}(t_j;\z_0^{(k)}) - \z(t_j;\z_0^{(k)}) \right\Vert^2.
\end{equation}
Here $\z_{\text{NN}}^{\theta}(t_j;\z_0^{(k)})$ is computed by applying
a numerical integrator \cite{lambert1991numerical} (e.g., Runge–Kutta methods) to the equation $\dot{\z}(t) = F_{\text{NN}}(\z;\theta)$
starting at $\z(t_{j-1};\z_0^{(k)})$;
$\|\cdot\|$ is the standard Euclidean norm.
Typically, gradient-based optimization methods are used
to solve this minimization problem.

If no physics is involved, 
the above gives a general description of the pure data-driven approach.
If some principles of physics are known, 
the physics-informed data-driven approach
aims to embed the available physics into neural networks.
There are two typical ways of doing the embedding. 
One is to add regularization terms to the loss \eqref{eq:loss} that penalize $F_{\text{NN}}$
that do not follow the physics.
The other is to devise a network architecture so that $F_{\text{NN}}$ obeys 
the available physics for any $\theta$ and the optimal parameters are then found by
minimizing the same loss \eqref{eq:loss}.

\subsection{The GENERIC formalism}
The General Equation for Non-Equilibrium Reversible-Irreversible Coupling (GENERIC) formalism
provides a general mathematical framework 
describing \textit{beyond-equilibrium} thermodynamic systems
\cite{ottinger2005beyond} including both conservative and dissipative systems. 
As a consequence, any system described by Hamilton's equation or Poisson's 
equation can be written in the GENERIC formalism, as follows: 

\begin{equation} \label{def:GENERIC}
    \begin{split}
        &\dot{\z}(t) = L(\z)\frac{\partial E}{\partial \z}(\z) + M(\z)\frac{\partial S}{\partial \z}(\z), 
    \\
    \text{subject to} 
    \quad
    &L(\z)\frac{\partial S}{\partial \z}(\z)
    =  M(\z)\frac{\partial E}{\partial \z}(\z) = \bm{0},
    \\
    &L(\z) \text{ is skew-symmetric, i.e., } L(\z) = -L(\z)^\top \\
    &M(\z) \text{ is symmetric positive semi-definite.}
    \end{split}
\end{equation}
The term $L\frac{\partial E}{\partial \z}$ accounts for all the reversible (non-dissipative) phenomena of the system. 
In the classical mechanics, this term is equivalent to Hamilton and Poisson's equations of motion.
The operator $L$ is called the Poisson matrix 
and is required to be skew-symmetric.
The term $M\frac{\partial S}{\partial \z}$ accounts for the irreversible (dissipative) material properties of the system. This term was motivated by the Ginzburg-Landau equation, which can be used to describe critical dynamics of spatially extended systems. The operator $M$ is called the friction matrix
and is required to be symmetric positive semi-definite.
$E(\z)$ and $S(\z)$ are 
the system's total energy and entropy,
respectively. Under this framework, it can be checked that the energy of the system is conserved and the entropy of the system monotonically increases with respect to time, i.e.,
$\frac{dE(\z(t))}{dt} =0$ and $\frac{dS(\z(t))}{dt} \ge 0$,
corresponding to the first and second laws of thermodynamics, respectively.

\subsection{Deep Neural Networks}
We employ deep neural networks
as the basic components in our surrogate modellings
for $L, M, E, S$.
For simplicity of discussion, we shall focus on 
feed-forward neural networks throughout this work,
while any types of neural networks (e.g., ResNet) can easily be used in place of the feed-forward networks 
without difficulties.

A $L$-layer feed-forward neural network $f:\mathbb{R}^{d_\text{in}} \to \mathbb{R}^{d_\text{out}}$ is 
defined by $f(x) = z^{L+1}(x)$
where $z^{L+1}$ is constructed 
recursively according to
\begin{equation*}
z^{\ell}(x) = W^{\ell}\phi(z^{\ell-1}(x)) + b^{\ell}, \qquad 1 < \ell \le L+1,
\end{equation*}
starting with $z^1(x) = W^1x + b^1$.
Here, $W^\ell \in \mathbb{R}^{n_{\ell} \times n_{\ell-1}}$ is the weight matrix and $b^\ell \in \mathbb{R}^{n_{\ell}}$ is the bias vector in the $l$-th layer,
where we set $n_{0} = d_\text{in} = d$ and $n_{L} = d_\text{out}$.
$\phi(x)$ is a nonlinear activation function that is applied element-wise.
The activation function is assumed to have certain properties
so that the universal approximation theorem holds \cite{Cybenko_MoC89,siegel2020approximation,Mhaskar_Neurl96}. 
The collection of all weights and biases of the network
is denoted by 
$\theta = \{(W^1, b^1), \cdots, (W^L, b^L)\}.$
We note that 
a neural network can be a matrix-valued function by converting the output to
a matrix of proper size.

Let $f(x;\theta):\mathbb{R}^{d} \to \mathbb{R}$ be a $L$-layer neural network
and $g:\mathbb{R}^d\to \mathbb{R}^d$
be a differentiable function.
The gradient of $(f\circ g)(x):=f(g(x))$ with respect to $x$ is given by
\begin{equation} \label{NN:multiplicative}
        \nabla_{x} f(g(x);\theta) 
        = (\emph{\textbf{J}}g(x))^\top\left(W^LD_{L-1} \cdots W^2D_1W^1\right)^\top
        = (\emph{\textbf{J}}g(x))^\top \nabla_{x} f \circ g(x),
\end{equation}
where $D_j$ is a diagonal matrix whose $(i,i)$-entry is $\phi'(z^j_i(g(x)))$
for $1 \le j < L$ and $1 \le i \le n_j$,
and $\textbf{J}g(x)$ is the Jacobian of $g$ at $x$.
A key observation is that 
$\nabla_{x} f(g(x);\theta)$
belongs to the row space of $\textbf{J}g(x)$
due to the multiplicative structure. 
We refer to the equation \eqref{NN:multiplicative}
as the multiplicative structure of the gradient of neural networks.

\section{GENERIC Formalism Informed Neural Networks} \label{sec:methods}
Our goal is to design neural network architectures 
that satisfy the symmetry and degeneracy conditions of \eqref{def:GENERIC},
yet are sufficiently expressive to learn the underlying dynamics from data.
Also, we want the proposed neural networks to be easily adopted when some prior physics is available (in terms of the GENERIC).

\begin{figure}[ht]
    \centering
    \includegraphics[width=0.99\textwidth]{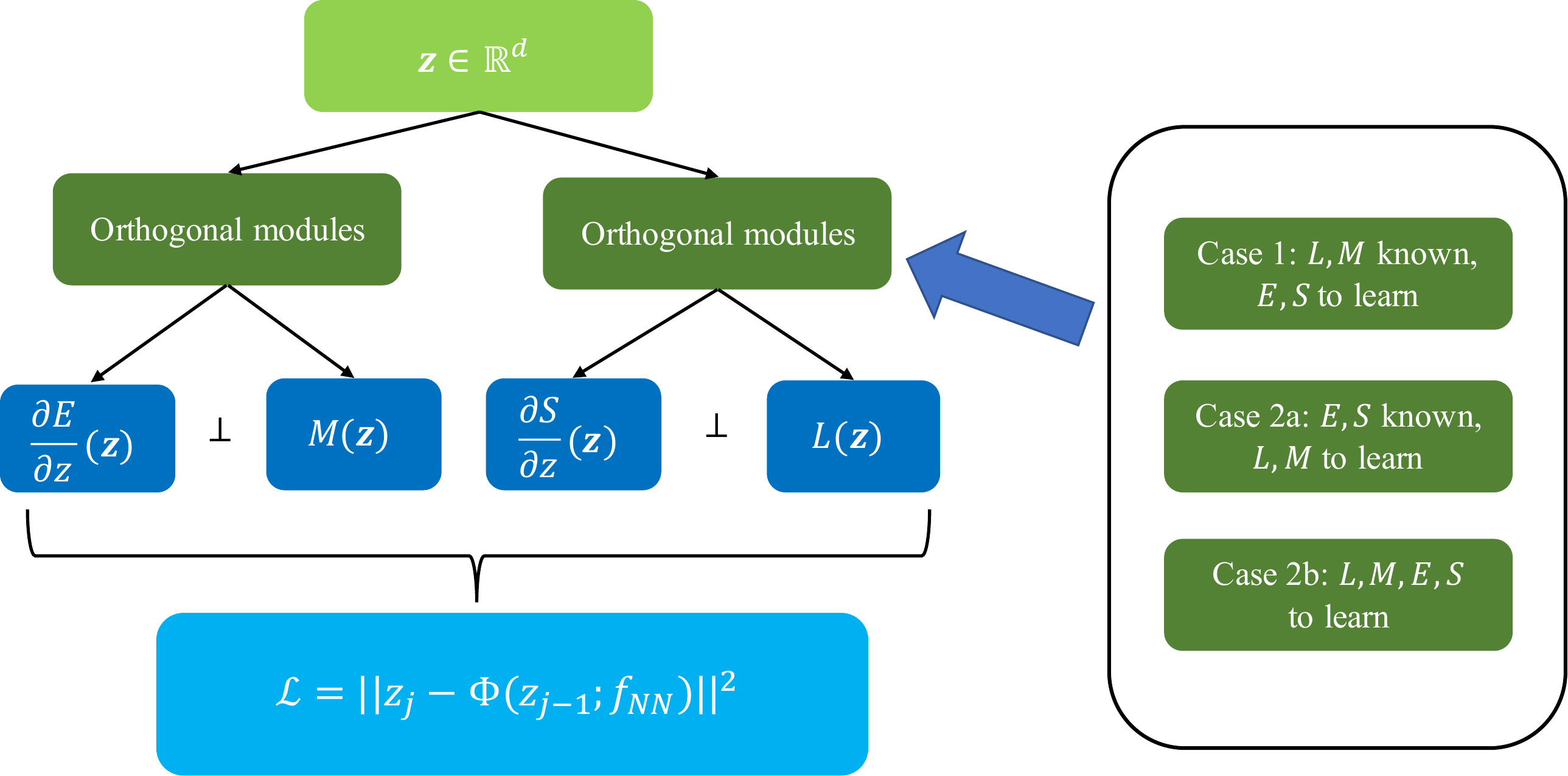}
    \caption{\textbf{Architecture of GFINNs.} GFINNs can be seen as a neural network framework, which automatically satisfies the laws of thermodynamics. The design of GFINNs follows the GENERIC formalism since we parameterize the (matrix-valued) functions $L,M,E,S$ as orthogonal neural network modules, which satisfy the consistency condition and structural condition in Eq.~\ref{def:GENERIC}. We propose different architectures, which can incorporate different physical information including the knowledge of $L,M$ or the knowledge of $E,S$.}
    \label{fig:architecture}
\end{figure}

Since the GENERIC formalism comprises two orthogonal modules,
each of which contains two components, 
we model each component using neural networks
that satisfy the required conditions.
The component-wise design not only allows the flexibility in incorporating prior physics (if any) but also results in a general framework of neural network modellings for the GENERIC formalism.
Here, prior physics implies a scenario where 
one or more
is known among $L$, $M$, $E$, $S$.
Although many possibilities can be discussed, 
for the sake of simplicity, 
we focus on the following scenarios:
\begin{itemize}
    \item Case 1: $L$ and $M$ are known.
    The goal is to approximate 
    $E$ and $S$.
    \item Case 2a: 
    $E$ and $S$ are unknown. 
    The goal is to approximate $L$ and $M$.
    \item Case 2b: 
    $L, M, E, S$ are unknown. 
    The goal is to approximate $L, M, E, S$.
\end{itemize}
All other scenarios can be handled without difficulties. 
The schematic of the proposed framework 
is shown in Figure~\ref{fig:architecture}.


Requiring only the conditions of \eqref{def:GENERIC} can be done quite easily,
while care is needed to ensure expressivity.
We illustrate this through the following example. 
Let 
\begin{align*}
    L(\z) = \begin{bmatrix} 0 & \ell(\z) & 0 & 0 \\ -\ell(\z) & 0 & 0 & 0\\ 0 & 0 & 0 & 0 \\ 0 & 0 & 0 & 0 \end{bmatrix},
    \qquad
    \ell(\z) \ne 0.
\end{align*}
Suppose the goal is to construct a proper neural network $S_{\text{NN}}$ satisfying $L(\z)\nabla S_{\text{NN}}(\z) = 0$.
Since $\text{null}(L(\z)) =  \text{span}\{e_3, e_4\}$, 
where $e_3 = (0, 0, 1,0)^\top$ and $e_4 = (0,0,0,1)^\top$, 
we have infinitely many choices for $S_{\text{NN}}$ that satisfy the condition.
For example, if $S_{\text{NN}}$ is any smooth function of $\z_3$ and does not depend on the other variables, 
the required condition is trivially satisfied. 
However, this class of functions does not capture the potential dependence of $\z_4$,
so that a sufficient expressiveness is not guaranteed. 
Similarly, if $S$ were given and the goal was to 
model a skew-symmetric matrix $L_{\text{NN}}$ that satisfies $L_{\text{NN}}(\z)\nabla S(\z) = 0$,
there are also infinitely many choices for $L_{\text{NN}}$ that satisfy the required condition, 
yet not all of them are sufficiently expressive.



For ease of discussion, in what follows, 
we denote a generic matrix function by $A$ 
representing either $L$ or $M$,
and a generic scalar function by $G$ representing either 
$S$ or $E$.
Also, the subscript ($G_{\text{NN}}$ or $A_{\text{NN}}$)
indicates the neural network model for the target quantity ($G$ or $A$).
Next, we discuss how to construct neural networks for 
modelling each component and present 
the corresponding universal approximation theorem 
under some assumptions.

\subsection{Case 1: \texorpdfstring{$L$}{L} and \texorpdfstring{$M$}{M} are known and \texorpdfstring{$E$}{E} and \texorpdfstring{$S$}{S} are unknown}
We consider the case where $L$ and $M$ are known, yet $E$ and $S$ are unknown.
Since the two generators are known, 
one might attempt to model $\nabla E$ and $\nabla S$
directly using neural networks following the pure data-driven approach \cite{ChenRBD18}.
However, 
since all vector functions are not gradient of a scalar function, 
we construct neural networks for $E$ and $S$
and compute their gradients by automatic differentiation that can be implemented by well-established programming packages e.g., Pytorch \cite{paszke2019pytorch} and Tensorflow \cite{abadi2016tensorflow}.
Furthermore, by construction, 
this approach 
allows one to not only predict the solution trajectories but also discover
physical quantities (energy and entropy)
from data.

Since we have multiple goals to achieve,
several challenges arise in developing 
neural network architectures with the desired properties. 
First of all, since we model $G$ (i.e., $E$ or $S$) instead of $\nabla G$,
we need to properly control the gradient of neural networks
so that the degeneracy condition of \eqref{def:GENERIC} holds. 
Motivated by the multiplicative structure \eqref{NN:multiplicative},
we introduce a tailored projection-like transformation $\mathcal{P}_A$ in
the very first layer of neural networks, 
which ensures the degeneracy condition
under some assumptions.
It is the transformation $\mathcal{P}_A$ 
that constitutes 
a core element in the network architecture
for $G$ (either $E$ and $S$)
assuming $A$ is known.

As a preparation for introducing the transformation $\mathcal{P}_A$ 
and also for the universal approximation theorem,
we make a couple of basic assumptions on $A$,
which will be justified in all the examples later.

\begin{assumption} \label{assumption1}
    Let $A(\z)$ be a $d\times d$ matrix-valued function defined on $\Omega \subset \mathbb{R}^d$.
    Let $q^j_{A}(\z)$, $j=1,\dots,n_A$,
    be an orthonormal basis of $\ker A(\z) := \{y \in \mathbb{R}^d : A(\z)y = 0\}$
    whose rank is $n_A$.
    We assume that 
    \begin{enumerate}
        \item $\ker A(\z)$ has constant rank
        $n_A \in \{1,\dots,d-1\}$
        in $\Omega$.
        \item 
        $\ker_{\text{inv}} A(\z):=\text{span}\{q_A^j:1\le j \le \tilde{n}_A\}$
        is 
        the largest subspace of $\ker A(\z)$
        whose rank is $\tilde{n}_A$
        which is also constant on $\Omega$
        such that  
        $q_A^j(\z)$, $j=1,\dots,\tilde{n}_A$,
        satisfy 
        \begin{equation} \label{eqn:assumption:cond1}
            \text{range}((\emph{\textbf{J}}q^j_{A}(\z))^\top) \subset \ker A(\z),
        \end{equation}
        where $\emph{\textbf{J}}q(\z)$ is the Jacobian matrix of $q(\z)$.
    \end{enumerate}
\end{assumption}

\begin{assumption} \label{assumption2}
    Let $A:\mathbb{R}^d \to \mathbb{R}^{d\times d}$ be a matrix-valued function
    satisfying Assumption~\ref{assumption1}
    and let $\hat{n}_A := n_A - \tilde{n}_A$.
    There exist
    real-valued differentiable functions $F^j_A(\z)$, $j=1,\dots,\hat{n}_{A}$,  on $\Omega$,
    satisfying 
    \begin{equation}
        \text{span}\{\nabla F_A^j(\z): j=1,\dots,\hat{n}_{A}\}
        \bigoplus
        \ker_{\text{inv}} A(\z)
        = \ker A(\z).
    \end{equation}
\end{assumption}

The degeneracy conditions of \eqref{def:GENERIC},
if it is interpreted in terms of linear algebra,
means that 
$\nabla G$ belongs to the kernel of $A$.
Since $\ker A(\z)$ depends on $\z$,
Assumption~\ref{assumption1}(i)
basically allows us to work on the fixed number of basis. 
This is a typical assumption made in order 
to make analysis go through, which is also used in \cite{golub1973differentiation} dated back in 1970s. 
Assumption~\ref{assumption1}(ii) and~\ref{assumption2} 
decompose $\ker A(\z)$
into two subspaces, one of which satisfies 
the invariance under differentiation 
in the sense of \eqref{eqn:assumption:cond1}.
The other subspace from Assumption~\ref{assumption2}
provides the key information on 
how to keep $\nabla G_{\text{NN}}$ in the kernel of $A(\z)$
as a function of $\z$.
The existence of the functions $F_A^j$
plays a key role in defining the transformation
$\mathcal{P}_A$.

With these assumptions, we define a transformation $\mathcal{P}_A$ as follows.
From Assumption~\ref{assumption1},
let $\tilde{Q}_A(\z) = [q_A^1(\z),\dots,q_A^{\tilde{n}_A}(\z)]$.
From Assumption~\ref{assumption2},
let $F_A(\z) = [F^1_A(\z), \dots, F^{\hat{n}_A}_A(\z)]^\top$.
Define the transformation $\mathcal{P}_A:\mathbb{R}^d \to \mathbb{R}^{n_A}$ by
\begin{equation} \label{def:transformation}
    \mathcal{P}_{A}(\z) = [\tilde{Q}_A^\top(\z)\z; {F}_{A}(\z)] \in \mathbb{R}^{n_A}.
\end{equation}
The first $\tilde{n}_A$ components of $\mathcal{P}_A(\z)$
is the orthogonal projection coefficients of $\z$
onto $\ker_{\text{inv}}A(\z)$
and the remaining components are  
the functions from Assumption~\ref{assumption2}.
The output of $\mathcal{P}_A$
is $n_A$-dimensional vector resulting in 
a dimension reduction from $d$ to $n_A$.

We are now in a position to present our 
neural networks for $G(\z)$ (either energy $E(\z)$ or entropy $S(\z)$).
For neural networks
$f(z;\theta):\mathbb{R}^d \to \mathbb{R}$,
we define 
\begin{equation} \label{def:ES-case1}
    \begin{split}
        G_{\text{NN}}(\z;\theta_A) = f(\mathcal{P}_A(\z); \theta_A).
    \end{split}
\end{equation}
It then follows from 
the multiplicative structure \eqref{NN:multiplicative}
and the properties of the transformation \eqref{def:transformation} stemmed from 
Assumptions~\ref{assumption1} and~\ref{assumption2}
that 
any $G_{\text{NN}}$ of the form \eqref{def:ES-case1}
satisfies 
$A(\z)\nabla G(\z) = \bm{0}$ 
for all $\z \in \Omega$.
We note that in general, 
the functions from Assumption~\ref{assumption2}
may not be readily available.
However, in Propositions~\ref{prop:gas} and~\ref{prop:dp},
we show that $F_A$ can be identified by extracting relevant components from the basis of $\ker A(\z)$ and then applying indefinite integration.

The remaining goal is to show 
the expressivity of the proposed neural network architecture \eqref{def:ES-case1}.
In order to show the universal approximation theorem,
an appropriate function class should be chosen 
in the first place
on which we show the universality.
Since the target function $G$ satisfying the degeneracy condition of \eqref{def:GENERIC}
depends highly on the properties of the kernel of an operator $A$,
a general function class 
requires some detailed characterizations of $\ker A(\z)$.
We thus confine ourselves to 
the function class $\mathcal{F}_A$ 
characterized by the transformation operator $\mathcal{P}_A$
together with the multiplicative structure \eqref{NN:multiplicative}.

\begin{definition}
    Let $A:\mathbb{R}^d \to \mathbb{R}^{d\times d}$ be a matrix-valued function
    satisfying Assumptions~\ref{assumption1} and~\ref{assumption2}.
    Let $\mathcal{F}_A$
    be the collection of differentiable functions $G$ defined on $\Omega \subset \mathbb{R}^d$ 
    whose gradient satisfies 
    \begin{equation} \label{def:F_A}
    \nabla G(\z) =
    (\emph{\textbf{J}}\mathcal{P}_A(\z))^\top \bm{c}_A \circ \mathcal{P}_A(\z),
    \qquad \forall \z \in \Omega,
    \end{equation}
    for some continuous function $\bm{c}_A:\mathbb{R}^{n_A} \to \mathbb{R}^{n_A}$.
\end{definition}


The function class $\mathcal{F}_A$ of Definition~\ref{def:F_A}
depends crucially on $A$.
For example, 
as shown in Corollary~\ref{cor:FA},
if $A$ is constant,
$\mathcal{F}_A$ contains 
all the continuously differentiable functions $G$
that satisfy $\nabla G \in \ker A(\z)$.

\begin{corollary} \label{cor:FA}
    Suppose $A$ is constant. 
    Then, 
    $\mathcal{F}_A$ is 
    the class $C^1(\Omega)$
    that consists of all differentiable functions whose gradient lies in $\ker A$
    and continuous on $\Omega$.
\end{corollary}
\begin{proof}
    Since $A$ is constant,
    so is $Q$.
    Also  
    $n:=n_A = \tilde{n}_A$, $\hat{n}_A = 0$,
    $\mathcal{P}_A(\z) = Q^\top \z$
    and $(\textbf{J}\mathcal{P}_A(\z))^\top= Q$.
    Let $\mathfrak{G}(\xi) := G(Q\xi)$
    where $\xi \in \mathbb{R}^n$.
    It then follows from 
    $\nabla G(\z) = Q\nabla_{\xi} G(Q\xi)\big|_{\xi = Q^\top \z}$
    that 
    $\nabla G(\z) 
    = Q\nabla_{\xi} \mathfrak{G}(\xi)\big|_{\xi=Q^\top \z}
    = Q\nabla_{\xi} \mathfrak{G} \circ \mathcal{P}_A(\z)$.
    By letting $\bm{c}_A(\xi) = \nabla_{\xi} \mathfrak{G}(\xi)$,
    the proof is completed. 
\end{proof}

With the target function class being defined,
we now show that 
the proposed neural network $G_{\text{NN}}$
defined in \eqref{def:ES-case1} 
is universal for 
the function class $\mathcal{F}_A$.

\begin{theorem}
    Suppose $A(\z)$ is a $d\times d$-matrix-valued function on $\Omega \subset \mathbb{R}^d$
    satisfying Assumptions~\ref{assumption1}
    and~\ref{assumption2}.
    For any $G(\z) \in \mathcal{F}_{A}$
    defined as in \eqref{def:F_A}
    and a small $\epsilon > 0$,
    there exists
    a neural network model $G_{\text{NN}}(\z)$
    defined as in \eqref{def:ES-case1}
    such that  
    \begin{equation*}
        \sup_{\z \in \Omega} 
        \|\nabla G(\z) - \nabla G_{\text{NN}}(\z) \| < \epsilon,
    \end{equation*}
    and $A(\z)\frac{\partial G_{\text{NN}}}{\partial \z}(\z)= \bm{0}$
   for all $\z \in \Omega$.
\end{theorem}
\begin{proof}
Note that 
$\frac{\partial G_{\text{NN}}}{\partial \z}(\z) = (\textbf{J}\mathcal{P}_A(\z))^\top \nabla f \circ \mathcal{P}_A(\z)$.
Since $G^* \in \mathcal{F}_A$, 
there exists a continuous function 
$\bm{c}_A^*$ such that 
the gradient of $G^*$ is expressed as 
$\nabla G^*(\z) = (\textbf{J}\mathcal{P}_A(\z))^\top \bm{c}_A^* \circ \mathcal{P}_A(\z)$.
Since $\textbf{J}\mathcal{P}_A(\z)$ is full rank on $\Omega$,
by invoking the universal approximation theorem of neural networks (e.g., \cite{li1996simultaneous,siegel2020approximation}), 
for any sufficiently small $\epsilon > 0$,
there exists a neural network $f$
satisfying 
$
\sup_{\z \in \Omega} \left\|
\nabla f \circ \mathcal{P}_A(\z) -\bm{c}_A^* \circ \mathcal{P}_A(\z) \right\| < \epsilon
$,
which completes the proof.
\end{proof}


An explicit neural network architecture in terms of width and depth may be given 
from existing works,
however, we simply rely on the well-established 
universal approximation theorem for neural networks 
\cite{Cybenko_MoC89,Mhaskar_Neurl96,li1996simultaneous,siegel2020approximation}.

The considered function class 
$\mathcal{F}_A$ from Definition~\ref{def:F_A}
may be too restricted to cover 
more general functions.
However, in all the examples of
Section~\ref{sec:experiments},
we show that all the assumptions hold and the underlying energy $E$ and entropy $S$ functions belong to
the function classes 
$\mathcal{F}_M$ and $\mathcal{F}_L$,
respectively.

\subsection{Case 2: \texorpdfstring{$L$}{L}, \texorpdfstring{$M$}{M}, \texorpdfstring{$E$}{E} and \texorpdfstring{$S$}{S} are unknown}
We consider only Case 2b 
where all the quantities $(L,M,E,S)$ are unknown,
since Case 2a is easily handled 
by letting $G_{\text{NN}} := G$.

In Case 2, 
a scalar function $G$ is modelled by a standard neural network unless it is known a priori. 
It then suffices to construct a matrix-valued neural network $A_{\text{NN}}$ from $G_{\text{NN}}$ that satisfies both the symmetry
and the degeneracy conditions.
Unlike Case 1, we do not need to control the gradient of $G_{\text{NN}}$ to be in $\ker A(\z)$.
Rather, we design $A_{\text{NN}}$ 
to satisfy $\nabla G_{\text{NN}} \in \ker A_{\text{NN}}(\z)$.
It turns out that 
one can easily adopt this property into 
neural network architectures
by exploiting skew-symmetric matrices. 

\begin{lemma} \label{lemma-Q}
    For $j=1,\dots,K$, let $S_j$ be a skew-symmetric matrix of size $d\times d$.
    For a differentiable scalar function $g(x):\mathbb{R}^d \to \mathbb{R}$,
    let $Q_g(x) \in \mathbb{R}^{K\times d}$ be a matrix-valued function whose $j$-th row is 
    defined to be 
    $(S_j\nabla g(x))^\top$.
    Then, $Q_g(x) \nabla g(x)= 0$ for all $x \in \mathbb{R}^d$.
\end{lemma}
\begin{proof}
    The proof readily follows from the fact that for any $y \in \mathbb{R}^d$, $y^\top S_j y = 0$
    since $S_j$ is skew-symmetric.
\end{proof}

For a neural network $G_{\text{NN}}$,
let $Q_{G_{\text{NN}}}(\z)$
be the matrix function defined through 
$K$ skew-symmetric matrices as in Lemma~\ref{lemma-Q},
which are trainable parameters.
Motivated by the spectral decomposition of 
either skew-symmetric or symmetric positive semi-definite
matrix,
we propose to model $A$ by 
\begin{equation} \label{def:L_M_NN}
    \begin{split}
        A_{\text{NN}}(\z) := (Q_{G_{\text{NN}}}(\z))^\top B_{\text{NN}}^A(\z)
        Q_{G_{\text{NN}}}(\z),
    \end{split}
\end{equation}
where $B_{\text{NN}}^A(\z)$ is skew-symmetric if $A= L$ 
and is symmetric positive semi-definite if $A = M$
that is modelled by another neural network 
different from $G_{\text{NN}}$.
In particular, 
we use two triangular matrix-valued neural networks $T_L(\z)$ and $T_M(\z)$ and set 
\begin{equation} \label{def:T_M}
    B_{\text{NN}}^L(\z) := (T_L(\z))^\top - T_L(\z),
    \qquad
    B_{\text{NN}}^M(\z) := (T_M(\z))^\top T_M(\z).
\end{equation}
By construction,
the symmetry and degeneracy conditions of \eqref{def:GENERIC} are automatically satisfied.

We note that 
if $A$ is either skew-symmetric or symmetric positive semi-definite, 
the spectral decomposition reads 
$A = Q^\top \Lambda Q$,
where $Q$ is orthogonal and $\Lambda$
is either skew-symmetric or diagonal.
The network architecture of \eqref{def:L_M_NN}
has a similar structure of that of 
the spectral decomposition,
yet, 
neither $Q_{G_{\text{NN}}}$ is orthogonal
nor $B_{\text{NN}}^A$ is the eigenvalue matrix.
Since $Q_{G_{\text{NN}}}$ is a matrix of size $K\times d$,
the rank of $A_{\text{NN}}$ is at most $\min\{K,d\}$.
Hence, $K$ is assumed to be greater than or equal to 
the rank of $A$.

Owing to the universal approximation theorem 
\cite{Cybenko_MoC89,Mhaskar_Neurl96}
of neural networks,
we show that 
the proposed neural network $A_{\text{NN}}$ of \eqref{def:L_M_NN} 
is sufficiently expressive enough to approximate the underlying target function $A$
under some mild conditions.

\begin{theorem}
    Suppose $\nabla G_{\text{NN}} := \nabla G$ is continuous on a compact set $\Omega \subset \mathbb{R}^d$,
    and a component of $\nabla G(\z)$ has nonzero values in $\Omega$.
    Let $A(\z)$ be
    either a skew-symmetric or symmetric positive semi-definite matrix-valued continuous function satisfying
    $A(\z)\nabla G(\z) = 0$ 
    for all $z \in \Omega$.
    For any $\epsilon > 0$, 
    there exists a neural network model $A_{\text{NN}}$ 
    of the form \eqref{def:L_M_NN}
    such that 
    $\sup_{\z \in \Omega} \|A - A_{\text{NN}}\| < \epsilon$
    and $A_{\text{NN}}(\z)\nabla G(\z) = \bm{0}$ in $\Omega$.
\end{theorem}
\begin{proof}
    Without loss of generality, let $(\nabla G(\z))_1 \ne 0$ for all $z \in \Omega$.
    Since $A\nabla G = 0$, 
    the column space $\text{col} A(\z)$ is spanned by at most $d-1$ independent basis.
    For $k=1,\dots,d-1$, let $P_k$ be a skew-symmetric matrix of size $d$ such that 
    $(P_k)_{1,k+1} = 1$, $(P_k)_{k+1,1}=-1$ and $(P_k)_{ij}=0$ otherwise.
    Let $\tilde{q}_k(\z) = P_k\nabla S(\z)$. 
    \begin{claim}
        $\tilde{q}_k(\z)$, $1\le k <d$, are linearly independent and 
        $$
        \text{col}(A(\z)) \subseteq \text{span}\{\tilde{q}_k(\z) : 1\le k < d\}.
        $$
    \end{claim}
    \begin{proof}[Proof of Claim]
        Assuming $\sum_{k=1}^{d-1} c_k \tilde{q}_k(\z) = 0$,
        it suffices to show $c_k = 0$ for all $k$. 
        Observe that 
        \begin{align*}
        \sum_{k=1}^{d-1} c_k \tilde{q}_k(\z) = \begin{bmatrix}
        \sum_{k=1}^{d-1}c_k (\nabla G(\z))_{k+1} &
        -c_1(\nabla G(\z))_1 & \cdots & -c_{k-1}(\nabla G(\z))_1 
        \end{bmatrix}^\top.
    \end{align*}
        Since $(\nabla G(\z))_1 \ne 0$, $c_1=\cdots = c_{k-1} = 0$, 
        which proves the linearly independence.
        Note also that since $P_k$ is skew-symmetric,
        $\langle \tilde{q}_k(\z), \nabla G(\z)\rangle = 0$ for all $k$.
        The second claim is followed from the relationship
        \begin{equation*}
            \text{col}(A(\z)) \subseteq \left(\text{span}\{\nabla G(\z)\}\right)^\perp = \text{span}\{\tilde{q}_k(\z) : 1\le k < d\},
        \end{equation*}
        where the equality holds because the two spaces have the same rank.
    \end{proof} 
    
    Note that $A$ is either 
    skew-symmetric 
    or symmetric positive semi-definite
    of rank $r$ with $2 \le r < d$. 
    We shall consider only the case when
    $A$ is skew-symmetric. 
    The other case can be proved similarly. 
    Thus, it can be decomposed as 
    $A(\z) = Q(\z)\Lambda(\z)Q^\top(\z)$,
    where $Q(\z)$ is a matrix of size $d\times r$ such that $Q^\top(\z) Q(\z) = I \in \mathbb{R}^{r\times r}$
    and  $\Lambda(\z)$ is skew-symmetric of size $r\times r$ defined by
    \begin{equation*}
        (\Lambda(\z))_{ij} = \begin{cases}
        \lambda_i(\z) & \text{if }  j = i+1, \\
        -\lambda_i(\z) & \text{if } j = i-1, \\
        0 & \text{otherwise},
        \end{cases}
    \end{equation*}
    where $\lambda_i(\z) > 0$.
    Thus, there exists a matrix $R(\z) \in \mathbb{R}^{(d-1)\times r}$ such that 
    $Q(\z) = \tilde{Q}(\z)R(\z)$,
    where $\tilde{Q}(\z) = [\tilde{q}_1(\z), \cdots, \tilde{q}_{d-1}(\z)]$.
    Since $Q^\top(\z) \tilde{Q}(\z)$ is full rank and $r \le d-1$,
    a solution to $Q^\top(\z) \tilde{Q}(\z) R(\z) = I$ always exists.
    
    Since $A(\z)$ is continuous on $\Omega \subset \mathbb{R}^d$,
    so does $\tilde{\Lambda}(\z):=R(\z)\Lambda(\z)R^\top(\z)$. 
    From the universal approximation theorem of 
    neural networks, 
    there exists
    a skew-symmetric matrix-valued neural network
    $B_\text{NN}^A(\z)$ such that 
    \begin{equation} \label{eqn:universal-approximation}
        \|(B_\text{NN}^A(\cdot))_{ij} - (\tilde{\Lambda}(\cdot))_{ij}\|_{C^0(\Omega)} < \frac{\epsilon}{C(d-1)r},
        \qquad \forall 1\le i < j < d,
    \end{equation}
    where $C = \sup_{\z \in \Omega} \|\tilde{Q}(\z)\|^2$.
    Note that $C$ is a constant that depends only on $\nabla G(\z)$ and $\Omega$.
    Therefore,
    \begin{align*}
        \|A(\z) - A_{\text{NN}}(\z)\|
        &= 
        \|Q(\z)\Lambda(\z)Q^\top(\z) - \tilde{Q}(\z)B_\text{NN}^A(\z)\tilde{Q}^\top(\z)\| \\
        &=
        \|\tilde{Q}(\z)R(\z)\Lambda(\z)R^\top(\z)\tilde{Q}^\top(\z) - \tilde{Q}(\z)B_\text{NN}^A(\z)\tilde{Q}^\top(\z)\| \\
        &\le \|\tilde{Q}(\z)\|^2 \|R(\z)\Lambda(\z)R^\top(\z) - B_\text{NN}^A(\z)\|.
    \end{align*}
    Since $\|\cdot\| \le \|\cdot\|_{F}$,
    where $\|\cdot\|_F$ is the Frobenius norm,
    it follows from \eqref{eqn:universal-approximation} that 
    $\|A(\z) - A_{\text{NN}}(\z)\| \le \epsilon$ for all $\z \in \Omega$, 
    which completes the proof.
\end{proof}



Assuming Case 2b, 
the GENERIC formalism informed neural networks (GFINNs) comprise of four neural networks:
$E_{\text{NN}}, S_{\text{NN}}, B_{\text{NN}}^L, B_{\text{NN}}^M$
together with sets of trainable parameters forming 
skew-symmetric matrices used in $Q_{E_{\text{NN}}}$
and $Q_{S_{\text{NN}}}$.
Apart from parameters for the four neural networks,
the number of parameters from the skew-symmetric matrices 
is $2Kd(d-1)$.
When the dimension $d$ is large, 
the number of trainable parameters grows $\mathcal{O}(Kd^2)$,
which may cause some computational challenges.
If this is the case, 
sparse parameterization can be applied to reduce the number of parameters. 
For example, 
each skew-symmetric matrix could be sparsely parameterized by only $s$ nonzero parameters,
which reduces the number of trainable parameters from $\mathcal{O}(Kd^2)$ to $\mathcal{O}(Ks)$.
A similar sparse parameterization can be applied in
modelling $B_{\text{NN}}^A$ as well.



In summary, 
GFINNs
consist of four components, 
$E_{\text{NN}}$, $S_{\text{NN}}$,
$L_{\text{NN}}$, $M_{\text{NN}}$
that form two modules; E-M module ($E_{\text{NN}}$, $M_{\text{NN}}$)
and L-S module ($L_{\text{NN}}$,$S_{\text{NN}}$).
When $A$ is known, 
we model $G$ using a neural network $G_{\text{NN}}$
of the form \eqref{def:ES-case1}.
When $G$ is known,
we model $A$ using a neural network $A_{\text{NN}}$
of the form \eqref{def:L_M_NN}.
When no physics is known,
we model $G$ using a standard neural network
$G_{\text{NN}}$
and then model $A$ using a neural network $A_{\text{NN}}$
of the form \eqref{def:L_M_NN}
with $G_{\text{NN}}$.
As a result, we obtain 
a general framework of designing neural networks 
for the GENERIC formalism 
with the flexibility of incorporating available physics,
thanks to the component-wise network modelling.

\section{Numerical Examples} \label{sec:experiments}

We demonstrate the performance of GFINNs
on three benchmark problems.
To compare against other methods, 
we also report the results obtained by
GNODEs \cite{lee2021machine},
SPNNs \cite{hernandez2021structure},
and SDENets \cite{dietrich2021learning}.
Implementation details can be found in Appendix \ref{sec:implementation}.
We note that 
SPNNs
and GNODEs 
can only be applied in case 1
and case 2, respectively,
while GFINNs cover all the cases.

\textbf{Nonuniqueness.}
By the multiplicative structure,
the GENERIC formalism allows multiple modules 
resulting in the same dynamics. 
That is, 
given $(E,S,L,M)$,
there are infinitely many 
$(\Tilde{E},\Tilde{S},\Tilde{L},\Tilde{M})$ satisfying
$\Tilde{L}(\z)\frac{\partial \Tilde{E}}{\partial \z}(\z) + \Tilde{M}(\z)\frac{\partial \Tilde{S}}{\partial \z}(\z)
    =
    L(\z)\frac{\partial E}{\partial \z}(\z) + M(\z)\frac{\partial S}{\partial \z}(\z)$.
As a matter of fact,  
if we let $\Tilde{G}$ be an affine transformation of $G$
with a slope coefficient $a$,
by letting $\Tilde{A} := a^{-1}A$,
we obtain new modules resulting in the same dynamics.
Because of the nonuniqueness, 
the inferred quantities $E_{\text{NN}}$ and $S_{\text{NN}}$
may look different from the target quantities.
We therefore calibrate the inferred quantities 
to make them look similar to the ground truth values.
The calibration is done by finding aforementioned affine transformations using some target values.
We apply the calibration only for the visualization purpose
to demonstrate the discovery of the physical quantities by GFINNs.

\textbf{GENERIC formalism under fluctuations.} 
Fluctuations can be included in the GENERIC formalism \cite{grmela1997dynamics,ottinger2005beyond},
resulting in a stochastic differential equation (SDE)
of the form
\begin{equation} \label{def:GENERIC_stoc}
    \begin{split}
    &d\z  
    = \mu(\z)dt + \sigma(\z)d\bm{W}_t,
    \\
    \text{subject to} 
    \quad
    &L(\z)\frac{\partial S}{\partial \z}(\z)
    =  \sigma(\z)\frac{\partial E}{\partial \z}(\z) = \bm{0},
    \\
    &
    L(\z) \text{ is skew-symmetric, i.e., } L(\z) = -L(\z)^\top, \\
    &
    \sigma(\z)\sigma(\z)^T = 2k_B M(\z),
    \end{split}
\end{equation}
where 
$\mu(\z) := L(\z)\frac{\partial E}{\partial \z}(\z) + M(\z)\frac{\partial S}{\partial \z}(\z)+k_B\frac{\partial}{\partial \z}\cdot M$,
$\bm{W}_t$ is a multicomponent Wiener process
and 
$k_B$ is the Bolzmann constant, which controls the magnitude of fluctuation. 
When the fluctuations are eliminated by letting $k_B\to 0$,
we recover 
\eqref{def:GENERIC} from \eqref{def:GENERIC_stoc}. 
Here $\frac{\partial}{\partial \z} \cdot M$ represents the divergence of $M$ as a tensor field, i.e.,
$\frac{\partial}{\partial \z} \cdot M = \sum_{i=1}^d\sum_{k=1}^d\frac{\partial M_{ik}}{\partial{\z_k}}\bm{e}_i$
where 
$\bm{e}_i$'s are the standard basis vectors 
in $\mathbb{R}^d$.
The consistency condition of \eqref{def:GENERIC_stoc}
implies the conservation of energy and fluctuation-dissipation theorem.

In the stochastic setting, the goal is to infer the drift and diffusion terms $\mu, \sigma$ 
from observed data. 
Let $Z(t,\omega;Z_0)$
be the solution to the SDE \eqref{def:GENERIC_stoc},
where $\omega$ denotes that $Z(t,\omega;Z_0)$ is a random variable and possesses the initial condition $Z(t_0,\omega;Z_0) = Z_0$
with probability one.
The data are then referred to 
$N_{\text{traj}}$ sample paths
of the solution to \eqref{def:GENERIC_stoc},
each of which has 
different initial states, $\{Z_0^{(k)}\}$, 
sampled from a probability distribution.
The $k$-th sample path is then written as 
$\{t_j, \z_j^{(k)}\}_{j=0}^{T}$
where 
$\z_j^{(k)}:=Z(t_j,\omega^{(k)};Z_0^{(k)})$
and $\omega^{(k)}$ is a realization of outcome.

GFINNs for the SDE \eqref{def:GENERIC_stoc} consist of 
the components $A_{\text{NN}}, G_{\text{NN}}, Q_{E_{\text{NN}}}$, $T_M$ from \eqref{def:L_M_NN} and \eqref{def:T_M},
which further construct 
$\mu_{\text{NN}}$ and  $\sigma_{\text{NN}}$
as follows:
\begin{equation} \label{def:GFINN-stoc}
    \begin{split}
        \mu_{\text{NN}}(\z) &:= L_{\text{NN}}(\z)\frac{\partial E_{\text{NN}}}{\partial \z}(\z) + M_{\text{NN}}(\z)\frac{\partial S_{\text{NN}}}{\partial \z}(\z) + k_B\frac{\partial}{\partial \z}\cdot M_{\text{NN}}(\z),
        \\
        \sigma_{\text{NN}}(\z) &:=\sqrt{2k_B}(T_M(\z)Q_{E_{\text{NN}}}(\z))^\top.
    \end{split}
\end{equation}
Both $\mu_{\text{NN}}(\z)$ and $\sigma_{\text{NN}}$ are naturally defined thanks to the spectral structure of $A_\text{NN}$. 
Since each component obeys the required conditions of \eqref{def:GENERIC},
GFINNs \eqref{def:GFINN-stoc} satisfy the consistency conditions of \eqref{def:GENERIC_stoc}.

\textbf{Loss function.} 
In the deterministic examples, the loss function is set to the mean squared error (MSE) defined in \eqref{eq:loss} 
together with the Runge-Kutta second/third order integrator \cite{lambert1991numerical}.

In the stochastic example, the loss function is set to the negative log-likelihood function \eqref{eq:loss_stoc}
together with the Euler-Maruyama integrator \cite{kloeden1992stochastic}. 
The same loss function is 
also used in \cite{dietrich2021learning, schneider2014maximum},
which is defined by
\begin{equation}\label{eq:loss_stoc}
\begin{split}
    \mathcal{L}(\theta) = -\frac{1}{N_{\text{traj}}} \sum_{k=1}^{N_{\text{traj}}} \frac{1}{T} \sum_{j=1}^{T} 
    \log p(\z_{j}^{(k)}|\z_{j-1}^{(k)},\theta),
\end{split}
\end{equation}
where 
$p(\z_{j}^{(k)}|\z_{j-1}^{(k)},\theta)$
is the probability density function of 
multivariate normal distribution with 
mean $\Delta t_j\mu_{\text{NN}}(\z_{j-1}^{(k)})$ 
and covariance matrix $2k_B\Delta t_j M_{\text{NN}}(\z_{j-1}^{(k)})$
evaluated at $\z_{j}^{(k)}-\z_{j-1}^{(k)}$.
Here $\Delta t_j := t_j - t_{j-1}$.

\textbf{Evaluation metric.}
The performance quality of learned dynamics 
is measured by 
a closeness    
between 
unseen ground truth trajectories 
that are not used in training,
and 
trajectories of inferred dynamics.
This is often referred to  
as generalization or test error.

Let $N_{\text{test}}$ be 
the number of unseen test trajectories.
For $k = 1,\dots,N_{\text{test}}$,
let $\bm{Z}^{(k)} \in\mathbb{R}^{T\times d}$ 
be the matrix representing 
the $k$-th \textit{unseen} trajectory
of ground-truth, 
whose $j$-th row, denoted by $\bm{Z}^{(k)}_j$,
is the state at time $t_j$.
Similarly,
let
$\bm{\Tilde{Z}}^{(k)}
\in\mathbb{R}^{T\times d}$ 
be the $k$-trajectory
matrix of
learned dynamics
whose initial state
is the same as the one of $\bm{Z}^{(k)}$.

In the deterministic case,
the metric we use for closeness is the mean squared error (MSE):
\begin{equation} \label{def:MSE-test}
    \text{MSE}(t_j)
    =\frac{1}{N_{\text{test}}}\sum_{k=1}^{N_{\text{test}}}
    \frac{1}{d}
    \sum_{l=1}^d (\bm{Z}_{jl}^{(k)} - \bm{\tilde{Z}}_{jl}^{(k)})^2.
\end{equation}

In the stochastic case, the metric we use for closeness is the squared sliced Wasserstein-2 distance \cite{deshpande2018generative},
which is defined through random projections.
Let $\{\bm{u}_m\}_{m=1}^M$ be a set of 
$M$ vectors randomly uniformly sampled from the unit hypersphere $\mathbb{S}^{d-1}$,
where we set $M = 100$ for implementation. 
For each $j$ and $m$, 
$\langle \bm{Z}_{j}^{(k)},\bm{u}_m \rangle$ 
and $\langle \bm{\Tilde{Z}}_{j}^{(k)}, \bm{u}_m \rangle$ 
are always assumed to be sorted with respect to the index $k$.
The squared sliced Wasserstein-2 distance (SW) is then defined by 
\begin{equation*}
    \text{SW}(t_j) =
    \frac{1}{N_{\text{test}}}
    \sum_{k=1}^{N_{\text{test}}}
    \frac{1}{M}
    \sum_{m=1}^M
    |\langle 
    \bm{Z}_{j}^{(k)} - \bm{\Tilde{Z}}_{j}^{(k)},
    \bm{u}_m
    \rangle|^2.
\end{equation*}

\subsection{Two gas containers exchanging heat and volume}
We consider the gas container example from 
\cite{shang2020structure}. 
Two gas containers are allowed to exchange heat and volume with a wall in the middle. The state variable
is $\z = (q,p,S_1,S_2) \in \mathbb{R}^4$, 
where $q$, $p$ represent the position and momentum of the moving wall, and 
$S_1$, $S_2$ represent the entropy of the gases in two containers. 
The energy of the whole system is  $E(\z) = \frac{p^2}{2m} + E_1 + E_2$,
where $E_i = (\frac{e^\frac{S_i}{Nk_B}}{\hat{c}V_i})^{\frac{2}{3}}$,
$V_1 = q$, $V_2 = 2-q$ and $\hat{c}= \left(\frac{4\pi m}{3h^2N}\right)^{\frac{3}{2}}\frac{e^{\frac{5}{2}}}{N}$,
which follows from the Sackur–Tetrode equation \cite{schroeder1999introduction} for ideal gases. 
$m$ is the mass of the wall, $N$ is the number of gas particles, $h$ is the Planck constant and $k_B$ is the Boltzmann constant. We fix the units such that $m = Nk_B = \hat{c} = 1$. The entropy of the system is $S(\z) = S_1 + S_2$.
The evolution equation is described by a system of  ordinary differential equations (ODEs):
\begin{equation} \label{def:LM-gc}
\begin{split}
    &\begin{pmatrix} \dot{q} \\ \dot{p} \\ \dot{S_1} \\ \dot{S_2} \end{pmatrix}=\begin{pmatrix} \frac{p}{m} \\ \frac{2}{3}(\frac{E_1}{q} - \frac{E_2}{2 - q}) \\
    \frac{\alpha}{T_1}(\frac{1}{T_1}-\frac{1}{T_2}) \\
    -\frac{\alpha}{T_2}(\frac{1}{T_1}-\frac{1}{T_2})\end{pmatrix}=L\frac{\partial E(\z)}{\partial \z} + M(\z)\frac{\partial S(\z)}{\partial \z},
\end{split}
\end{equation}
where 
$L = \begin{bmatrix}
\bm{S} & \bm{0} \\ \bm{0} & \bm{0}
\end{bmatrix}
$,
$
M = \begin{bmatrix}
\bm{0} & \bm{0} \\ \bm{0} & \bm{T}
\end{bmatrix}
$, 
$\bm{S} = \begin{bmatrix}
0 & 1 \\ -1 & 0
\end{bmatrix}
$, 
$\bm{T} =  \alpha \begin{bmatrix}
T_1^{-2} & -(T_1T_2)^{-1} \\ -(T_1T_{2})^{-1} & T_2^{-2}
\end{bmatrix}$, 
$T_i = \frac{\partial E_i}{\partial S_i}$, 
$\bm{0}$ is the zero matrix of size $2\times 2$
and 
$\alpha$ is the parameter determining the strength of heat exchange which we set to 10.

The following proposition
shows that 
the governing equation \eqref{def:LM-gc} 
of the gas container problem 
satisfies all the assumptions of Section~\ref{sec:methods}.
\begin{proposition} \label{prop:gas}
    Let $M(\z)$ be the matrix 
    defined in \eqref{def:LM-gc}.
    Let 
    \begin{align*}
        F_M(\z) = e^{CS_1}q^{-\frac{2}{3}}
        +
        e^{CS_2}(2-q)^{-\frac{2}{3}},
    \end{align*}
    where $C = \frac{2}{3Nk_B}$.
    Then, $\ker M(\z) = \text{span}\{e_1, e_2, \nabla F_M(\z)\}$, where
    \begin{equation*}
        \nabla F_M(\z) =
        (T_3, 0, T_1, T_2)^\top,
        \qquad
        T_3 = 
        -\frac{2}{3}e^{\frac{2}{3}CS_1}q^{-\frac{5}{3}}
        +\frac{2}{3}e^{\frac{2}{3}CS_2}(2-q)^{-\frac{5}{3}}.
    \end{equation*}
    Let $\hat{q}(\z) = (0, 0, \bar{T}_1, \bar{T}_2)^\top$
    where 
    $\bar{T}_i = \frac{T_i}{\sqrt{T_1^2+T_2^2}}$.
    Then, 
    $\{e_1,e_2,\hat{q}(\z)\}$ is an orthonormal basis for $\ker M(\z)$
    and $\text{range}((\emph{\textbf{J}}\hat{q}(\z))^\top) \not\subset \ker M(\z)$.
    Furthermore, $E \in \mathcal{F}_M$ as 
    \begin{equation*}
        \nabla E(\z) = (\emph{\textbf{J}}\mathcal{P}_M(\z))^\top \bm{c}_{M} \circ \mathcal{P}_M(\z),
    \end{equation*}
    where 
    $\mathcal{P}_M(\z) = (q, p, F_M(\z))^\top$, 
    $\bm{c}_M(\xi_1,\xi_2,\xi_3) = (1, \xi_2/m, 1)^\top$,
    and 
    $(\emph{\textbf{J}}\mathcal{P}_M(\z))^\top = [e_1, e_2, \nabla F_M(\z)]$.
\end{proposition}
\begin{proof}
    The proof directly follows from a straight forward calculation.
\end{proof}

We use 80 trajectories starting from 
$t_0 = 0$ to $t_T = 8$ with $\Delta t_j := t_j - t_{j-1} = 0.02$ ($T=400$) as training data
and use another 20 trajectories as test data. 
The initial conditions of both training and testing trajectories are uniformly sampled from $[0.2,1.8]\times[-1,1]\times[1,3]\times[1,3]$.
Since neural networks may be 
sensitive to how parameters are initialized,
we run ten independent simulations 
and report some ensembles out of it.
As an effort to make a fair comparison,
neural networks used for each method
have a roughly similar number of parameters.
The detailed architectures are summarized in
Table~\ref{tab:model_params}
of Appendix~\ref{sec:implementation}.
\begin{figure}[ht]
    \centering
    \includegraphics[width=0.99\textwidth]{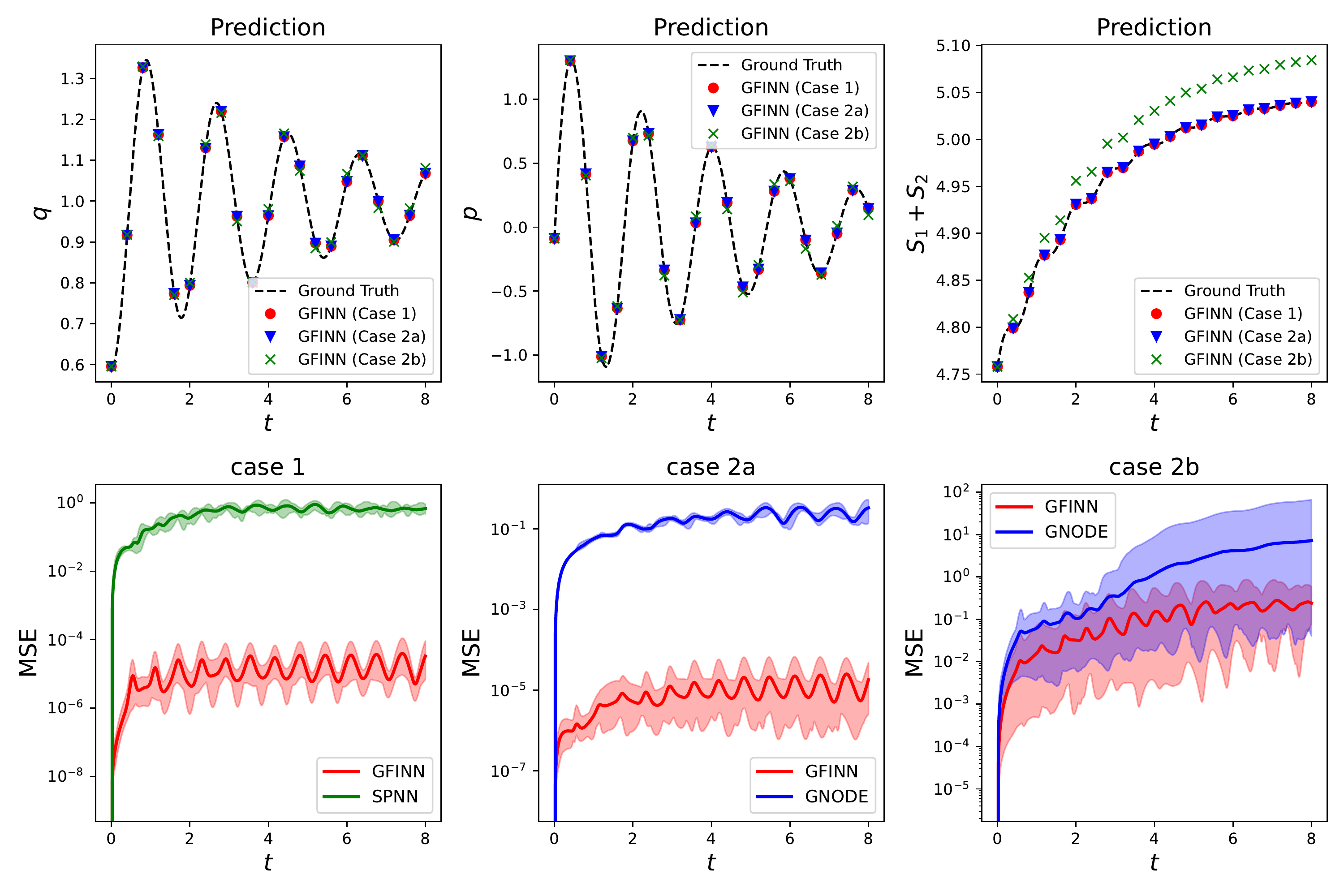}
    \caption{\textbf{Results for the gas container example.} (\textbf{Top-left and Top-middle}) Predicted position $q$ and momentum $p$ for the wall. (\textbf{Top-right}) Predicted entropy $S_1 + S_2$ of the air in two gas containers. (\textbf{Bottom-left}) The prediction MSE of GFINN and SPNN in case 1. (\textbf{Bottom-middle and bottom-right}) The prediction MSE of GFINN and GNODE in case 2a and case 2b. The results are obtained by taking the mean of 10 simulations, while the upper and lower boundary for shaded region represents the highest and lowest error in 10 simulations.
    }
    \label{fig:gc_error}
\end{figure}

In the top row of Figure~\ref{fig:gc_error}, 
we plot 
one of 20 test trajectories
as dashed-lines,
together with 
the corresponding trajectory of 
the learned dynamics by GFINNs 
as symbols;
the circle ($\circ$), 
the inverted triangle ($\triangledown$)
and the cross ($\times$) marks 
correspond to 
case 1, case 2a and case 2b, respectively. 
Here, GFINNs are the one with the smallest 
MSE summed over time among 10 simulations.
Since the state variable lies in $\mathbb{R}^4$,
the $q$, $p$ and $S_1+S_2$ trajectories are reported.
We clearly see that 
in case 1 and case 2a,
the trajectories of GFINNs
are indistinguishable to 
the ground truth trajectory,
while in case 2b, 
they start to deviate from the truth trajectory
as time increases.
This is expected 
as some underlying physics in terms of 
the GENERIC formalism is known in case 1 and case 2a, while no physics is known in case 2b.

To compare against other methods,
in the bottom row of Figure~\ref{fig:gc_error},
we report the means of the MSE \eqref{def:MSE-test}
of GFINNs, GNODEs and SPNNs 
with respect to time.
Each shaded region 
represents the range between the maximum and the minimum of the MSEs from ten simulations.
We clearly observe that 
the mean of the MSEs by GFINNs 
is much lower than 
those by SPNNs and GNODEs in all the cases.
In both case 1 and case 2a, 
all the 10 MSEs of GFINNs are 
significantly lower than those of the other comparisons.
This again demonstrates that by incorporating physical knowledge into neural networks, 
GFINNs can achieve much higher predictive accuracy.
In case 2b, 
the mean MSE of GFINNs is at least 
one order magnitude smaller than
the one of GNODEs 
in almost all times.
These results indicate 
that to achieve 
good performance,
it is not enough to just 
enforce the GENERIC conditions,
but a sufficient expressivity is
also required 
to capture the underlying dynamics.



\begin{figure}[ht]
    \centering
    \includegraphics[width=0.99\textwidth]{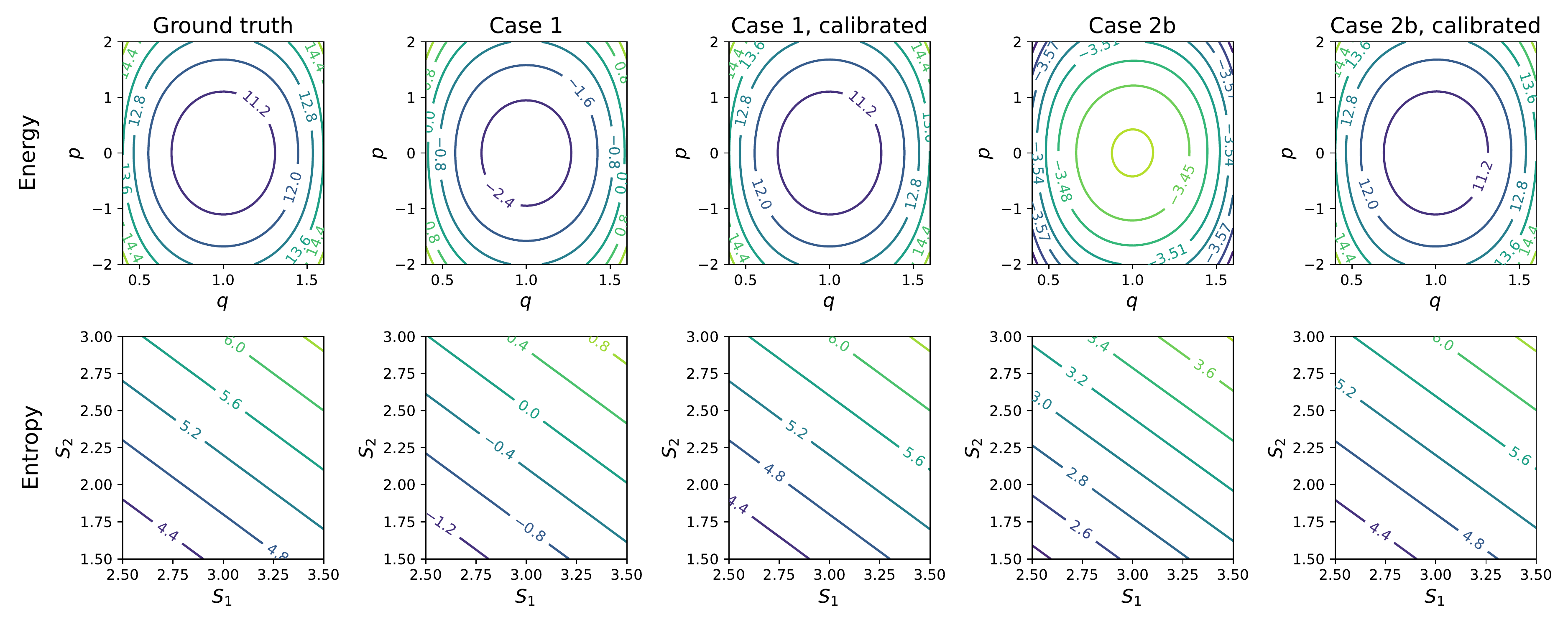}
    \caption{\textbf{Learned energy $E$ and entropy $S$ of the gas container example by GFINN.} In both case 1 and case 2b, the GFINNs can recover the correct energy and entropy of the system when the undetermined scaling and translation factor is calibrated with the ground truth. We show the energy surface at the hyperplane $\{z:S_1=2.5,S_2=2.5\}$ as a function of $p,q$ and the entropy surface at $\{z:p=0,q=1\}$ as a function of $S_1$ and $S_2$.}
    \label{fig:gc_contour}
\end{figure}

In Figure~\ref{fig:gc_contour},
we plot 
the contours of energy and entropy
functions from 
both the ground truth
and GFINNs in all the cases.
Due to the nonuniquness of 
the GENERIC formalism, 
a proper calibration is applied.
We see that the calibrated contours of 
both the energy and entropy
by GFINNs are indistinguishable
to those by the ground truth.
This demonstrates the discovery of the energy and entropy by GFINNs from data.

\subsection{Thermoelastic double pendulum}
We consider the two-dimensional finite thermoelastic double pendulum example \cite{romero2009thermodynamically}. The state variable is 
$\z = (\bm{q}_1, \bm{q}_2, \bm{p}_1, \bm{p}_2, S_1, S_2) \in \mathbb{R}^{10}$, 
where $\bm{q}_i, \bm{p}_i \in \mathbb{R}^2$ represent the position, momentum of the $i$-th mass, and
$S_i$ represents the entropy of the $i$-th spring. We denote the length of two springs as 
$\lambda_1 = \Vert\bm{q}_1\Vert$, $\lambda_2 = \Vert\bm{q}_2 - \bm{q}_1\Vert$. 
The total energy of the system is $E(\z) = \frac{\Vert\bm{p}_1\Vert^2}{2} + \frac{\Vert\bm{p}_2\Vert^2}{2} + E_1 + E_2$, where $E_1$ and $E_2$ are the internal energy of both springs, defined by
$E_i = \frac{1}{2} (\log\lambda_i)^2 + \log\lambda_i + e^{S_i-\log\lambda_i}-1$.
The entropy of the system is $S(\z) = S_1 + S_2$. 
The evolution equation for the problem is 
then given by a system of ODEs:
\begin{equation} \label{eqn-dp}
\begin{split}
    \begin{pmatrix} \dot{\bm{q}_1} \\\dot{\bm{q}_2} \\ \dot{\bm{p}_1} \\ \dot{\bm{p}_2} \\ \dot{S_1} \\ \dot{S_2} \end{pmatrix}
    &=\begin{pmatrix} \bm{p}_1 \\ \bm{p}_2 \\
    -\frac{\partial}{\partial \bm{q}_1}(E_1 + E_2) \\
    -\frac{\partial}{\partial \bm{q}_2}(E_1 + E_2) \\
    T_1^{-1}T_2 - 1 \\ T_1T_2^{-1} - 1
    \end{pmatrix}=L\frac{\partial E(\z)}{\partial \z} + M(\z)\frac{\partial S(\z)}{\partial \z},
\end{split}
\end{equation}
where
$T_i = \frac{\partial E_i}{\partial S_i}$, 
$\bm{0}$ and $\bm{1}$
are $2\times 2$ matrices whose elements
are all 0 and 1, respectively, 
$\bm{0}_{m_1 \times m_2}$
is a matrix of size $m_1\times m_2$
whose elements are all 0, 
\begin{equation} \label{def:LM-DP}
    L =\begin{pmatrix} 
    \bm{0}_{4\times 4}
    & \bm{S}
    \\
    -\bm{S}^\top & \bm{0}_{6\times 6}
    \end{pmatrix}, 
    \quad 
    M(\z) = \begin{pmatrix}
    \bm{0}_{8\times 8} &  \bm{0}_{8\times 2}
    \\
    \bm{0}_{2\times 8} & \bm{T}(\z)
    \end{pmatrix}, \quad 
    \bm{S} = \begin{pmatrix}
    \bm{1} & \bm{0} & \bm{0} \\
    \bm{0} & \bm{1} & \bm{0} 
    \end{pmatrix}, \quad 
    \bm{T} = \begin{pmatrix}
    \frac{T_2}{T_1} & -1 \\
     -1 & \frac{T_1}{T_2}
    \end{pmatrix}.
\end{equation}
The governing equation \eqref{eqn-dp} for the problem
again satisfies all the assumptions in Section \ref{sec:methods},
which is shown in the following proposition.
\begin{proposition} \label{prop:dp}
    Let $M(\z)$ be the matrix defined 
    in \eqref{def:LM-DP}.
    Let $F_M(\z) = 
        -\|\bm{q}_1\|e^{S_1}
        -\|\bm{q}_2-\bm{q}_1\|e^{S_2}$.
    Then, $\ker M(\z) = \text{span}\{e_1,\dots,e_8,
    \nabla F_M(\z) \}$
    and 
    $\nabla F_M(\z) =
        (T_3; T_4; \bm{0}_{4\times 1}; T_1; T_2)$
    where 
    \begin{equation*}
        T_3 := -\frac{\bm{q}_1}{\|\bm{q}_1\|}e^{S_1}
        +\frac{\bm{q}_2-\bm{q}_1}{\|\bm{q}_2-\bm{q}_1\|}e^{S_2},
        \qquad
        T_4 := 
        -\frac{\bm{q}_2-\bm{q}_1}{\|\bm{q}_2-\bm{q}_1\|}e^{S_2}.
    \end{equation*}
    Let $\hat{q}(\z) = (\bm{0}_{8\times 1};  \bar{T}_1; \bar{T}_2)$
    where 
    $\bar{T}_i = \frac{T_i}{\sqrt{T_1^2+T_2^2}}$.
    Then, 
    $\{e_1,\dots,e_8,\hat{q}(\z)\}$ is an orthonormal basis for $\ker M(\z)$
    and $\text{range}((\emph{\textbf{J}}\hat{q}(\z))^\top) \not\subset \ker M(\z)$.
    Furthermore, $E \in \mathcal{F}_M$ 
    and $\nabla E(\z) = (\emph{\textbf{J}}\mathcal{P}_M(\z))^\top \bm{c}_{M} \circ \mathcal{P}_M(\z)$,
    where 
    $\mathcal{P}_M(\z) = (\bm{q}_1,\bm{q}_2,\bm{p}_1,\bm{p}_2, F_M(\z))^\top$, 
    $\bm{c}_M(\xi) = 
    (\bm{0}_{4\times 1}; \xi_5; \dots; \xi_8; 1)$
    and 
    $(\emph{\textbf{J}}\mathcal{P}_M(\z))^\top = [e_1, \dots, e_8, \nabla F_M(\z)]$.
\end{proposition}
\begin{proof}
    The proof directly follows from a straight forward calculation.
\end{proof}

 We generate 100 trajectories 
 from $t_0 = 0$ to $t_T = 40$ with $\Delta t = 0.1$
 ($T = 400$),
 whose initial conditions are sampled uniformly from $[0.9,1.1]\times[-0.1,0.1]\times[2.1,2.3]\times[-0.1,0.1]\times[-0.1,0.1]\times[1.9,2.1]\times[0.9,1.1]\times[-0.1,0.1]\times[0.9,1.1]\times[0.1,0.3]$.
 We use 80 of them for training and the remaining for testing.

\begin{figure}[ht]
    \centering
    \includegraphics[width=0.99\textwidth]{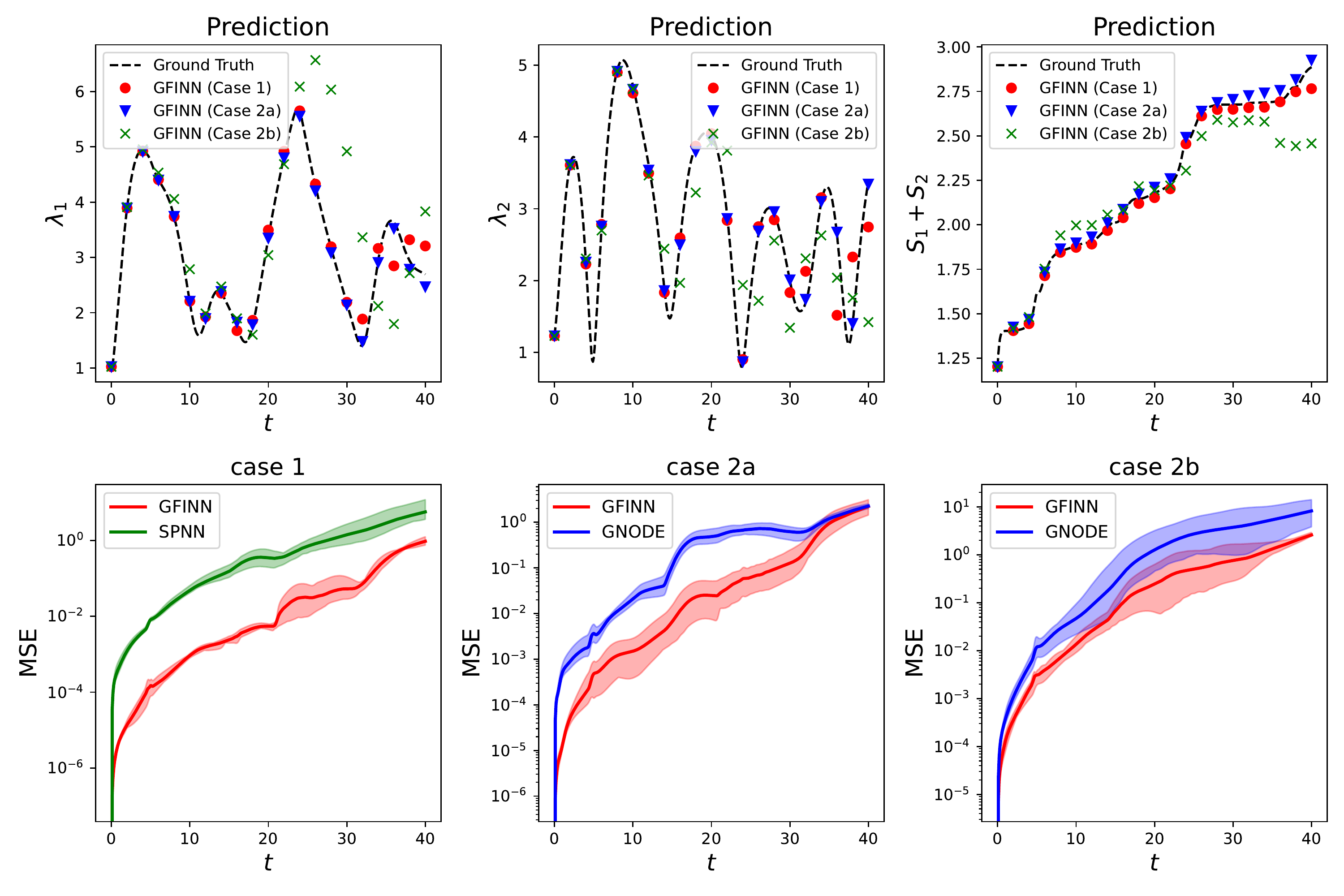}
    \caption{\textbf{Results for the thermoelastic double pendulum example.} (\textbf{Top}) Predicted length of springs $\lambda_1,\lambda_2$ and the predicted total entropy $S_1 + S_2$. Better predictive performance is achieved when there is additional physical knowledge (case 1 and case 2a) beyond the GENERIC formalism. (\textbf{Bottom}) The MSE of GFINNs, SPNNs and GNODEs. By incorporating hard constraints into the neural network, GFINNs can produce better predictive performances compared to SPNNs and GNODEs.}
    \label{fig:DP_error}
\end{figure}

In the top row of Fig~\ref{fig:DP_error}, we again plot one of 20 test trajectories, together with the predicted trajectories of GFINNs in the same way as in the previous example. Since the state variable lies in $\mathbb{R}^{10}$, we plot the predicted trajectories of $\lambda_1$, $\lambda_2$ and $S_1 + S_2$. We see that the predicted length and entropy of the springs by GFINNs match the ground truth in case 1 and case 2a. In case 2b, due to the chaotic nature of the problem, the predicted trajectory starts to deviate from the ground after certain time.

In the bottom row of Fig~\ref{fig:DP_error}, we report the mean MSE of GFINN, SPNNs and GNODEs from ten simulations.
The shaded areas indicate
the maximum and minimum of MSEs 
as in the gas container example. We clearly observe that in all cases GFINNs can produce better predictive performance compared to other baseline methods. In case 1, all the 10 MSEs of GFINNs are significantly smaller than those of SPNNs. In case 2a and case 2b, the mean MSE of GFINNs is approximately one order of magnitude smaller than GNODEs most of the time. The only exception is in the time window $t\in[30,40]$, when the MSEs for GFINNs and GNODEs become similar. 

\subsection{Langevin equation} 
We consider the diffusion of a particle described by the Langevin equation. The state variable vector is 
$z = (q,p,S_e) \in \mathbb{R}^3$, where $q$, $p$ are 
the position and momentum of the particle, respectively, and $S_e$ is the entropy of the surrounding environment. 
A simple form of energy and entropy is assumed:
$E = \frac{p^2}{2} + S_e$ and $S = S_e$. 
The Langevin equation is then given by \eqref{def:GENERIC_stoc} with 
\begin{equation} \label{def:LM-lg}
    L =\begin{pmatrix} 0 & 1  & 0 \\ -1 & 0 & 0 \\ 0 & 0 & 0 \end{pmatrix}, \quad
    M(\z) = \begin{pmatrix}  0 & 0 & 0 \\  0 & \frac{1}{2} & -\frac{p}{2} \\ 0 & -\frac{p}{2} & \frac{p^2}{2}\end{pmatrix}, \quad \sigma(\z) = 
    \begin{pmatrix}
    0\\ 1\\ -p
    \end{pmatrix}.
\end{equation}
Note that we choose units such that $k_B = 1$ for notational simplicity. 

\begin{figure}[ht]
    \centering
    \includegraphics[width=0.99\textwidth]{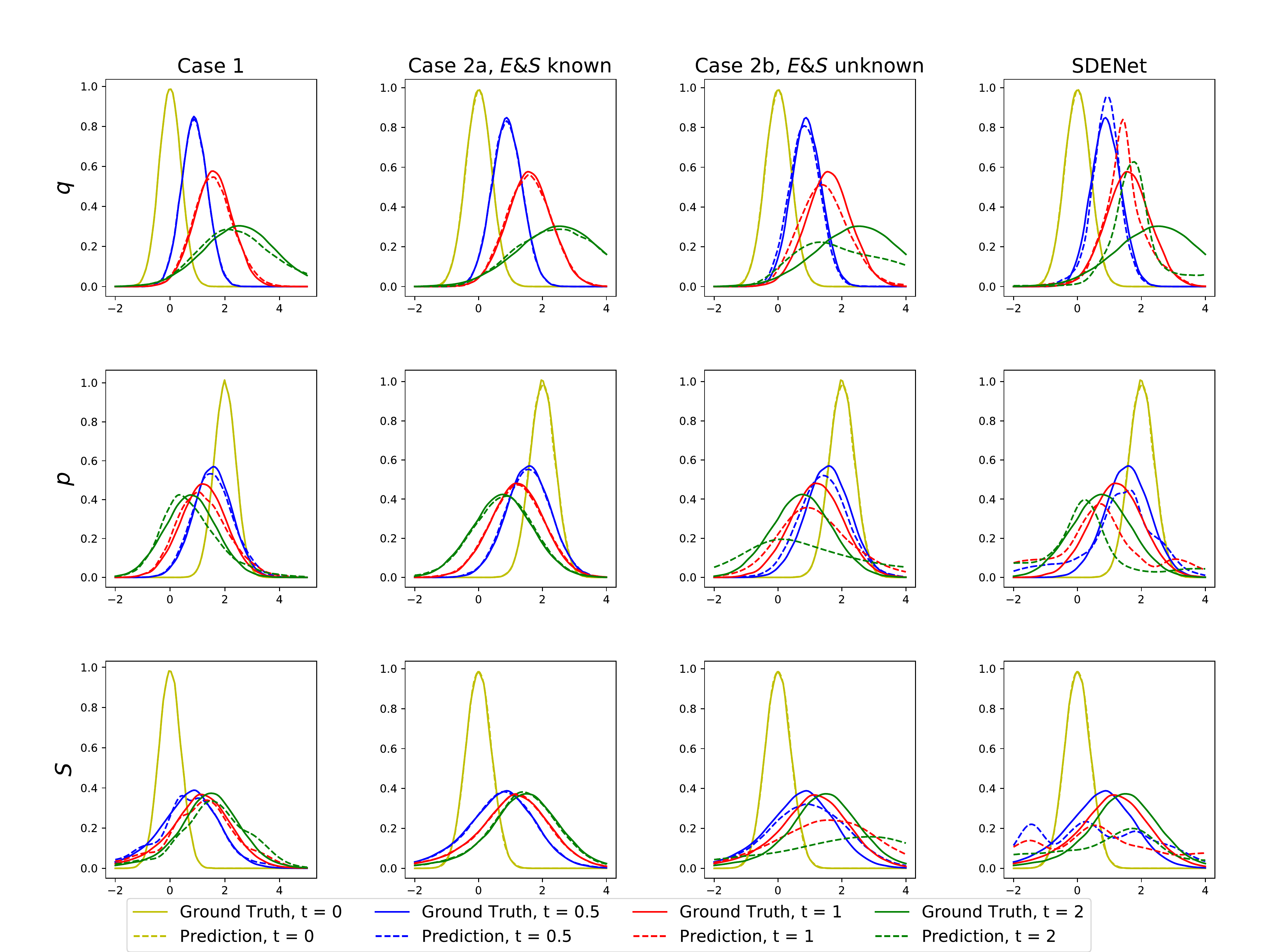}
    \caption{\textbf{Results for the Langevin equation example.} We plot the predicted and true probability distribution of $\z(t)$ at $t=0,0.5,1,2$ using both types of GFINNs. In case 2a, GFINNs can recover the correct distribution and give the correct prediction when extrapolating in time ($t = 2$). In case 1, there is a small discrepancy between the predicted and true distribution. For case 2b and SDENet, the prediction deviates from the ground truth after certain time.}
    \label{fig:LG_dist}
\end{figure}

We simulate 40 trajectories between $t_0 = 0$ to $t_T = 1$ with $\Delta t = 0.004$ for training. 
The initial state is 
randomly sampled from 
a normal distribution with 
mean $(0, 2, 0)^\top$
and covariance matrix $0.4^2I$,
where $I$ is the identity matrix of size $3$.
To evaluate the predictive performance of the model, we sample another 50,000 initial conditions from the same distribution, and solve \eqref{def:GENERIC_stoc} using both the ground truth $\mu, \sigma$ and inferred $\mu_{\text{NN}},\sigma_{\text{NN}}$ with the Euler-Maruyama method. Then we apply the Gaussian kernel density estimator 
    \cite{parzen1962estimation}
    to approximate the distribution of $\z$ and $\z_{\text{NN}}$ at time steps $t = 0, 0.5, 1, 2$, results of which can be found in Fig.~\ref{fig:LG_dist}. 
    Note that we infer dynamics from stochastic samples.
    In case 1 and case 2a, our model is not only able to recover the correct distribution of $\z(t)$, but also gives the correct prediction when extrapolating in time ($t = 2$), using a small amount of data. However, for case 2b and SDENet, even though a similar training error is achieved as shown in Table~\ref{tab:LG_loss}, the predicted trajectories do not match the ground truth, which indicates that prior physical knowledge is crucial to make the data-efficient inference. Still, the prediction error of GFINNs is smaller than that of SDENets in all cases
    as shown in Table~\ref{tab:LG_loss}.

\begin{table}[htbp]
\centering
\begin{tabular}{|c|c|c|c|c|}
\hline
 & Case 1 & Case 2a & Case 2b & SDENet \\ \hline
 Training loss & -2.7856 & -2.7844 & -2.6910 & -2.7690 \\ \hline
  Prediction error, $t = 0.5$ & $1.3\times 10^{-2}$ & $4.6\times 10^{-4}$ & $7.6\times 10^{-2}$ & $4.0\times 10^{-1}$ \\ \hline
 Prediction error, $t = 1$ & $6.7\times 10^{-2}$ & $7.4\times 10^{-4}$ & $7.4\times 10^{-1}$ & 1.1 \\ \hline
\end{tabular}
\caption{\textbf{Training loss (negative log-likelihood) and prediction error (squared sliced Wasserstein-2 distance) of GFINNs.} Even though the training losses of the three cases are similar, the prediction error shows significant differences, which reflects the different generalization capabilities of different models. In all cases, GFINNs produce a lower prediction error compared to SDENets.}
\label{tab:LG_loss}
\end{table}

\appendix
\section{Implementation of architecture} \label{sec:implementation}
\subsection{GNODEs}
The GNODEs \cite{lee2021machine} parameterize $L,M$ 
based on bracket structure and use standard neural networks for modelling $E, S$.
The GNODEs' architecture is designed for case 2b
and their architecture has a similarity with GFINNs in that 
two independent neural networks $E_{\text{NN}}, S_\text{NN}$ are used to parameterize $E$ and $S$.
However, a key difference lies in 
the network architectures for $L$ and $M$.
Specifically, GNODEs model $L$ and $M$ by
$L^{\text{GNODE}}$ and $M^{\text{GNODE}}$
whose $(\alpha,\beta)$ components are defined by
\begin{equation*}
    L_{\alpha\beta}^{\text{GNODE}}(\z) = \sum_{\gamma}\xi_{\alpha\beta\gamma}\nabla S_\text{NN}(\z)_{\gamma}, \hspace{0.2cm}
    M_{\alpha\beta}^{\text{GNODE}}(\z) =\sum_{\mu,\nu,m,n}\Lambda_{\alpha\beta}^mD_{mn}\Lambda_{\mu\nu}^n\nabla E_\text{NN}(\z)_{\beta}\nabla E_\text{NN}(\z)_{\nu}
\end{equation*}
where $\xi$ is a 3d skew-symmetric tensor, $\Lambda^m$ are 2d skew-symmetric matrices, $D$ is a positive semi-definite matrix. 
Due to the bracket structure, the GENERIC conditions are enforced under such parameterization. 
However, $L^{\text{GNODE}}$ and $M^{\text{GNODE}}$ are merely functions of $\nabla S_\text{NN}(\z)$ and $\nabla E_\text{NN}(\z)$ respectively, which can be seen as the underlying model assumption that differs from GFINNs. 
In examples where such assumption holds, GNODEs may achieve good performance
because they incorporate stronger physical prior knowledge. However, when $A$ depends on not only $\nabla G(\z)$ but also on some other functions of $\z$, GNODEs cannot represent the underlying governing equations due to the lack of expressivity. 
In \cite{lee2021machine}, the mean squared error \eqref{eq:loss} is also used as the loss function.

\subsection{SPNNs}
\cite{hernandez2021structure}
proposed SPNNs,
which parameterize the gradient of the energy and entropy
assuming $L$ and $M$ are known (case 1).
The loss function for SPNNs 
is defined as
\begin{equation*}
    \mathcal{L}(\theta) = \frac{1}{N_{\text{traj}}} \sum_{k=1}^{N_{\text{traj}}} \frac{1}{T} \sum_{j=1}^{T} 
    \left(\left\Vert\z_{\text{NN}}(t_j;\z_0^{(k)};\theta) - \z(t_j;\z_0^{(k)}) \right\Vert^2 + \lambda\left( \left\Vert L(dS)_{\text{NN}}\right\Vert^2 + \left\Vert M(dE)_{\text{NN}}\right\Vert^2\right)\right),
\end{equation*}
where 
$(dE)_{\text{NN}}$ and $(dS)_{\text{NN}}$ are standard feed-forward neural networks (FNNs) parameterizing $\nabla E$ and $\nabla S$,
$L(dS)_{\text{NN}}$ and $M(dE)_{\text{NN}}$ are both evaluated at $\z(t_j;\z_0^{(k)})$, $\lambda$ is the hyperparameter which controls the scale of the soft penalty,
and 
$\z_{\text{NN}}(t_j;\z_0^{(k)};\theta)$ is computed by applying
a numerical integrator to the equation $\dot{\z}(t) = F_{\text{NN}}(\z;\theta) = L(\z) (dE)_{\text{NN}}(\z) + M(\z) (dS)_{\text{NN}}(\z)$
starting at $\z(t_{j-1};\z_0^{(k)})$, 

However, the architecture we used for comparison is slightly different from that in the original paper, i.e., we model $E,S$ instead of $\nabla E$ and $\nabla S$ as neural networks. By parameterizing $E,S$ as neural networks, the surrogate model can learn a dynamical system more effectively, due to the inductive bias that $\nabla E$ and $\nabla S$ should satisfy the Clairaut's theorem. Another reason for our choice is that we parameterize $E,S$ as neural networks in GFINNs, and the comparison should be made fair by control of variables.

In this paper we choose $\lambda$ from $\{0.01,0.1,1\}$ and pick the ones that give the lowest MSEs in the two deterministic problems.
\subsection{SDENet}
The SDENets parameterize the drift and diffusion term $\mu$ and $\sigma$ as two independent FNNs and use the negative log-likelihood function \eqref{eq:loss_stoc} as the loss function.  

\begin{table}[htbp]
\centering
\begin{tabular}{|c|c|c|c|c|c|c|c|c|}
\hline
\multicolumn{2}{|c|}{\multirow{2}{*}{Problem}} & \multicolumn{2}{c|}{Gas container} & \multicolumn{2}{c|}{Double pendulum} & \multicolumn{3}{c|}{Langevin equation} \\ \cline{3-9}
\multicolumn{2}{|c|}{} & GFINN & GNODE & GFINN & GNODE & GFINN & \multicolumn{2}{c|}{SDENet} \\ \hline
\multirow{4}{*}{Layers}
 & $L$ & -/5/5 & - & -/5/5 & - & -/5/5 & \multirow{2}{*}{$\mu$} & \multirow{2}{*}{5}\\ \cline{2-7}
 & $M$ & -/5/5 & - & -/5/5  & -  & -/5/1 & &\\ \cline{2-9}
 & $E$ & 5/-/5 & 5 & 5/-/5 & 5 & 5/-/5 &\multirow{2}{*}{$\sigma$} & \multirow{2}{*}{5}\\ \cline{2-7}
 & $S$ & 5/-/1 & 1 & 5/-/5 & 5  & 5/-/1 & &\\ \hline
\multirow{4}{*}{Width}
 & $L$ & -/30/30 & - & -/30/30 & - & -/30/30 & \multirow{2}{*}{$\mu$} &\multirow{2}{*}{30}\\ \cline{2-7}
 & $M$ & -/30/30 & - & -/30/30 & - & -/30/- & &\\ \cline{2-9}
 & $E$ & 30/-/30 & 30 & 30/-/30 & 30 & 30/-/30 & \multirow{2}{*}{$\sigma$} &\multirow{2}{*}{30}\\ \cline{2-7}
 & $S$ & 30/-/- & - & 30/-/30  & 30  & 30/-/- & &\\ \hline
\end{tabular}
\caption{\textbf{Model architecture.} The number of layers and width of the neural network is slightly different across three cases considered in the paper (case 1, case 2a and case 2b). We list all of them in the table separated by the slash. The layers and width for each component of SPNNs is set to be the same as GFINNs. The hyperbolic tangent (tanh) activation functions are used for all the models. The optimizer is chosen to be Adam \cite{adam2015} with learning rate $0.001$. For the gas container and double pendulum examples, we train the models using mini-batch with batch size $100$ for $5\times 10^{5}$ iterations. For the Langevin equation example, we train the models using full batch for $5\times 10^{4}$ iterations. The important dimension $K$ of skew-symmetric matrices $Q_{S_{\text{NN}}}(\z)$ and $Q_{E_{\text{NN}}}(\z)$ in the parameterization of $L_{\text{NN}}$ and $M_{\text{NN}}$ are chosen to be 5 and 4 respectively in all examples. }
\label{tab:model_params}
\end{table}

\section*{Acknowledgments}
 We would like to acknowledghe support by DOE PhILMs (no. DE- SC0019453) and OSD/AFOSR MURI grant FA9550-
20-1-0358. We would like to acknowledge the helpful discussion with Dr. Xin Bian, Dr. Zhen Li and Dr. Chensen Lin. We would like to thank Dr. Kookjin Lee, Dr. Nat Trask and Dr. Panos Stinis for providing the codes of the GNODE paper.

\bibliographystyle{RS}
\bibliography{references}

\end{document}